\documentclass[letterpaper,11pt]{article}
\usepackage[usenames,dvipsnames,svgnames,table]{xcolor}
\usepackage{amsmath}
\usepackage{bbold}
\usepackage[all]{xy}
\usepackage{graphicx}
\usepackage{enumitem}
\usepackage{marginnote}
\usepackage{amsthm}
\input{amssym}
\usepackage{hyperref}
\hypersetup{
	colorlinks=true,
	linktoc=all,
	linkcolor=blue,
}
\usepackage[capitalize,nameinlink,noabbrev]{cleveref} 
\usepackage{cleveref}
\crefname{subsection}{subsection}{subsections}
\newtheorem{theorem}{Theorem}[section]
\newtheorem{lemma}[theorem]{Lemma}
\newtheorem{claim}{}[theorem]
\newtheorem{cor}[theorem]{Corollary}

\newtheorem{observation}[theorem]{Observation}
\theoremstyle{definition}
\newtheorem*{definition}{Definition}
\newtheorem*{notation}{Notation}
\newtheorem*{note}{Note}

\title{Splitter Theorems for Graph Immersions}
\author{Matt DeVos
\and
Mahdieh Malekian}
\date{}

\newcommand{\col}{black}

\begin{document}
\maketitle
\begin{abstract}
We establish splitter theorems for graph immersions for two families of graphs, $k$-edge-connected graphs, with $k$ even, and 3-edge-connected, internally 4-edge-connected graphs. As a corollary, we prove that every $3$-edge-connected, internally $4$-edge-connected graph on at least seven vertices that immerses $K_5$ also has $K_{3,3}$ as an immersion.
\end{abstract}

\vspace{0.7cm}

\section{Introduction}
\label{intro}

Throughout the paper, we use standard definitions and notation for graphs as in \cite{MR1367739}. Let $G$ be a graph with a certain connectivity. One natural question is that whether there is a way to ``reduce" $G$ while preserving the same connectivity, and possibly also the presence of a particular graph ``contained" in $G$. Broadly speaking, in answering such questions, two types of theorems arise. In the first type, \emph{chain theorems}, one tries to ``reduce" the graph down to some basic starting point, which is typically a particular small graph, or a small family of graphs. The other type of theorems are \emph{splitter theorems}. Here, there is the extra information that another graph $H$ is properly ``contained" in $G$ , and both have a certain connectivity. The idea is then to ``reduce" $G$ to a graph ``one step smaller", while preserving the connectivity, and the ``containment" of $H$.

The best known such results are the ones in the world where the connectivity concerned is vertex-connectivity, the``reduction" is an edge-contraction or edge-deletion, and the ``containment" relation is being isomorphic to a minor. In this realm, the first chain result is a well-known fact that if a graph $G$ is $2$-connected, then for every edge $e \in E(G)$, either $G\setminus e$ or $G/e$ is $2$-connected. The next result, due to Tutte, is a classical result of a reduction theorem of chain variety. Here, a \emph{wheel} is a graph formed by connecting a single vertex to all vertices of a cycle.

\begin{theorem}[Tutte \cite{MR0140094}]
Let $G$ be a simple $3$-connected graph. Then there exists $e \in E(G)$ such that either $G\setminus e$ or $G/e$ is simple and $3$-connected, unless $G$ is a wheel.
\end{theorem}

Another classical result of chain type is that every simple 3-connected graph other than $K_4$ has an edge whose contraction results in a $3$-connected graph, see \cite{MR1373655}. 
There is a wide body of literature sharpening these results and extending them to other connectivity, see, for instance \cite{MR2467821, MR1892444}.

Reduction theorems of splitter variety for graph minors started with a result for $2$-connected graphs which appears in a more general form in the work of Brylawski \cite{MR0309764}, and Seymour \cite{MR0439663}. This result asserts that if $G, H$ are $2$-connected graphs, and $H$ is a proper minor of $G$, then there is an edge $e \in E(G)$ such that $G\setminus e$ or $G/e$ is $2$-connected, and has $H$ as a minor. The more famous splitter theorem is Seymour's Splitter Theorem for $3$-connected graphs, which asserts:

\begin{theorem}[Seymour \cite{MR579077}]
Let $G$ and $H$ be $3$-connected simple graphs where $H$ is a proper minor of $G$ and $|E(H)|\ge 4$. Assume that if $H$ is a wheel, then $G$ has no larger wheel minor. Then $G$ has an edge $e$ such that either $G\setminus e$ or $G/e$ is simple and $3$-connected, and contains $H$ as a minor.
\end{theorem}

There is an extremely wide body of literature extending these results to other connectivity and the realm of matroids, and binary matroids, see, for example, \cite{MR2900811,costalonga2018splitter,MR3641803}.

In this paper, however, we are not concerned with vertex-connectivity and minors, but rather the world of edge-connectivity and a less-explored type of containment, \emph{immersions}. A pair of distinct edges $xy, yz$  or the 2-edge-path $xyz$ is said to \emph{split off at} $y$ if we delete these edges and add a new edge $xz$. We say a graph $G$ \emph{immerses} $H$, or \emph{has an $H$-immersion}, and write $G \succeq H$ or $H \preceq G$, if a subgraph of $G$ can be transformed to a graph isomorphic to $H$ through a series of splitting pairs of edges.\footnote{This is sometimes called \emph{weak immersion}.} If $G \ncong H$, we may write $G \succ H$ or $H \prec G$.  Also, we say a vertex $v \in V(G)$ of even degree is \emph{completely split} if $\deg(v)/2$ consecutive splits are performed at $v$, and then the resulting isolated vertex $v$ is deleted.

In the world of edge-connectivity, and immersions, there is a chain theorem due to Lov\'{a}sz (\cite{MR1265492}, Problem $6.53$, see also \cite{MR1373655}). 

\begin{theorem}[Lov\'{a}sz \cite{MR1265492}]
Let $G$ be a $2k$-edge-connected graph. Then by repeatedly applying complete split  and edge-deletion it can be reduced to a graph on two vertices, with $2k$ parallel edges between them.
\end{theorem}

This theorem was later generalized by a significant theorem of Mader that is key in our proofs, and will be stated in  \cref{kecon}. The goal of this paper is to establish two splitter theorems for immersions, the first of which is an analogue of the aforementioned result of Lov\'{a}sz.

\begin{theorem}
\label{sp-thm1intro}
Let $G$ and $H$ be $2k$-edge-connected loopless graphs where $G \succ H$. Then there exists an operation taking $G$ to $G'$ so that $G'$ is $2k$-edge-connected and $G' \succeq H$, where the operation is either
\begin{itemize}
\item deleting an edge, or
\item splitting at a vertex of degree at least $2k+2$, or
\item completely splitting a vertex of degree $2k$.
\end{itemize}
\end{theorem}

In comparison with graph minors, the literature on splitter theorems for graph immersions is extremely sparse. Indeed, we only know of two significant papers concerning this, namely \cite{MR2216473,MR2646129}, where Ding and Kanno have proved a handful of splitter theorems for immersion for cubic graphs, and $4$-regular graphs. In particular, they have shown the following (see \cite{MR2646129}, Theorem 9):

\begin{theorem}[Ding, Kanno \cite{MR2646129}]
\label{DingKanno}
Let $G$ and $H$ be $4$-edge-connected $4$-regular loopless graphs where $G \succ  H$. Then there exists a vertex whose complete split takes $G$ to $G'$ so that $G'$ is $4$-edge-connected $4$-regular, and $G' \succeq H$.
\end{theorem}

Our second theorem, stated below, is similar to the first one, but is for a different type of connectivity, and generalizes \cref{DingKanno}. Here, $G$ is said to be \emph{internally $k$-edge-connected} if every edge-cut containing less than $k$ edges is the set of edges incident with a single vertex. Also, in the statement of the theorem, $Q_3$ denotes the graph of the cube, and $K_2^3$ denotes the graph on two vertices with three parallel edges between them.

\begin{theorem}
\label{sp-thm2intro}
Let $G$ and $H$ be $3$-edge-connected and internally $4$-edge-connected loopless graphs, where $G \succ H$. Assume $|V(H)|\ge 2$ and $(G, H) \ncong (Q_3, K_4), (Q_3, K_2^3)$. Then there exists an operation taking $G$ to $G'$ such that $G'$ is $3$-edge-connected, internally $4$-edge-connected, and $G' \succeq H$, where the operation is either
\begin{itemize}
\item deleting an edge, or
\item splitting at a vertex of degree at least $4$,
\end{itemize}
each followed by iteratively deleting any loops, and suppressing vertices of degree $2$.
\end{theorem}

As an immediate corollary of the above theorem, we derive the following chain type result for the family of $3$-edge-connected, internally $4$-edge-connected graphs.
\begin{cor}
	Let $G$ be a $3$-edge-connected and internally $4$-edge-connected graph where $|V(G)| \ge 2$. If $G\ncong Q_3, K_2^3$, then one of the operations in the statement of \cref{sp-thm2intro} may be applied to $G$, where the resulting graph is $3$-edge-connected and internally $4$-edge-connected.
\end{cor}

In the world of graph minors, an immediate simple consequence of Seymour's Splitter Theorem, first observed by Wagner\cite{MR1513158}, is that every $3$-connected graph on at least six vertices containing $K_5$ as a minor, has a $K_{3,3}$-minor. This fact is then used to obtain a precise structural description of graphs with no $K_{3,3}$-minor. In parallel to this, and as an application of  \cref{sp-thm2intro}, we will establish the following analogue of this result for immersions. The result will be a step towards understanding graphs with no $K_5$-immersion.

\begin{cor}
\label{corintro}
Let $G$ be a $3$-edge-connected and internally $4$-edge-connected graph where $G\succ K_5$. Then
\begin{enumerate}
\item
if $|V(G)|= 6$, then $G\succeq K_{3,3}$ or $G\cong K_{2,2,2}$.
\item
If $|V(G)|\ge 7$, then $G\succ K_{3,3}$.
\end{enumerate}
\end{cor}

The rest of the paper is organized as follows: In  \cref{kecon}, we state the preliminary definitions and key tools, and prove  \cref{sp-thm1intro}.  \cref{in4econ} is dedicated to the family of $3$-edge-connected, internally $4$-edge-connected graphs, and includes the proof of \cref{sp-thm2intro} and \cref{corintro}.

\section{$k$-edge-connected graphs, $k$ even}
\label{kecon}

We will assume the graphs are undirected and finite, which may have loops or parallel edges. For $X\subset V(G)$, we use $\delta_G(X)$ to denote the edge-cut consisting of all edges of $G$ with exactly one endpoint in $X$, the number of which is called the \emph{size} of this edge-cut, and is denoted by $d_G(X)$. When $G$ is connected we refer to both $X$ and $X^c (=V(G)\setminus X)$ as \emph{sides} of the edge-cut $\delta (X)$. An edge-cut is called \emph{trivial} if at least one side of the cut consists of only one vertex. We call an edge-cut of size one a \emph{cut-edge}. Note that graph is $k$-edge-connected (internally $k$-edge-connected) if every edge-cut (non-trivial edge-cut) has size at least $k$. For distinct vertices $x, y \in V(G)$, we let $\lambda_G({x,y})$ denote the maximum size of a collection of pairwise edge-disjoint paths between $x$ and $y$. By a \emph{$t$-vertex}, we mean a vertex with degree $t$. For $X$ and $Y$ distinct nonempty subsets of $V(G)$, we denote the set of edges with one end in $X$ and the other end in $Y$ by $E_G(X, Y)$, and its size by $e_G(X, Y)$. Whenever the graph concerned is clear from the context, we may drop the subscript $G$. If either set is a singleton, say \( X = \{x\} \) and \( Y = \{y\} \), we may abbreviate \( E(X, Y) \), \( e(X, Y) \) to \( E(x, Y) \), \( e(x, Y) \), or to \( E(x, y) \), \( e(x, y) \) if both sets are singletons. 
 In \cref{intro} the notion of graph immersions was introduced. Equivalently, we could say $H$ is immersed in $G$ if there is an injective mapping $\phi:V(H) \rightarrow V(G)$ and a path $P_{uv}$ from $\phi(u)$ to $\phi(v)$  in $G$ corresponding to every edge $uv$ in $H$ so that $P_{uv}$ paths are pairwise edge-disjoint. In this case, a vertex in $\phi(V(H))$ is called a \emph{terminal} of $H$-immersion. \footnote{It is worth mentioning that if $G\succeq H$ and the collection of paths $P_{uv}$ are internally disjoint from $\phi(V(H))$, it is standard in the literature to say $G$ \emph{strongly immerses} $H$. It would be then in contrast with the notion of \emph{weak immersion}, where $P_{uv}$ paths are not necessarily internally disjoint from $\phi(V(H))$. However, we are only studying weak immersion and, for the sake of simplicity, refer to it as immersion.}

We proceed by listing a couple facts and theorems which will feature in our proofs. The first is the observation below, {\color \col {commonly known as the submodularity property}}.

\begin{observation}
\label{cuts-edge-count}
Let $G$ be a graph, and let $X$ and $Y$ be distinct nonempty subsets of $V(G)$. Then, by counting the edges contributing to the edge-cuts, we have
$$d(X\cap Y)+d(X\cup Y)+2e(X^c\cap Y, X\cap Y^c)=d(X)+d(Y).$$
Observe that it also implies the following inequality $$d(X\cap Y)+d(X\cup Y)\le d(X)+d(Y).$$
\end{observation}

Another frequently used fact is the classical theorem of Menger. A proof may be found, for instance, in \cite{MR1367739}.

\begin{theorem}[Menger]
\label{Menger}
Let $G$ be a graph, and let $x$ and $y$ be distinct vertices of $G$. Then $\lambda_G(x,y)$ equals the minimum size of an edge-cut of $G$ separating $x$ from $y$.
\end{theorem}

The next theorem is a strong form of Mader's splitting theorem that is a key ingredient in our proofs.

\begin{theorem}[Frank \cite{MR1172364}, see also Mader \cite{MR499117}]
	\label{Mader}
	Let $G$ be a graph. Assume that for $s\in V(G)$, we have $d(s) \ge 4$, and that $s$ is not incident with a cut-edge. Then there exist $\lfloor d(s)/2 \rfloor $ pairwise disjoint pairs $\{e_i, f_i\} \subseteq \delta(s)$ for $i=1, \ldots, \lfloor d(s)/2 \rfloor$, so that in the graph $G_i$ resulting from splitting $e_i, f_i$,  for any $x,y\in V(G_i) \setminus s$, we have $\lambda_G(x,y)=\lambda_{G_i}(x,y)$, $i=1, \ldots, {d(s)}/{2}$.
\end{theorem}

Throughout the rest of this section, we assume $k$ is even. The main result of this section is \cref{sp-thm1intro} which is {\color \col{repeated}} below for convenience.

\begin{theorem}
\label{sp-thm1}
Let $k \ge 4$ be an even number, and suppose $G$ and $H$ are $k$-edge-connected loopless graphs, where $G\succ  H$. Then there exists an operation taking $G$ to $G'$ so that $G'$ is  loopless and $k$-edge-connected and $G'\succeq   H$, where an operation is either
\begin{itemize}
\item deleting an edge, or
\item splitting at a vertex of degree at least $k+2$, or
\item completely splitting a $k$-vertex.
\end{itemize}
\end{theorem}

Note that in order to have a splitter theorem for the family of $k$-edge-connected graphs, we do need to embrace completely splitting a $k$-vertex as one of our operations, since as soon as we do a split at a $k$-vertex, the resulting graph will have a trivial $(k - 2)$-edge-cut.

\cref{sp-thm1} will be proved through a series of lemmas.   Although this theorem concerns loopless $k$-edge-connected graphs, our lemmas will have to allow for a wider family of graphs that have a special flaw.  This is motivated by the operation of completely splitting a $k$-vertex $u$.  We will consider this operation as a sequence of $\frac{k}{2}$  splits at $u$ done one at a time, and we shall need to investigate some of the graphs appearing partway  through this process.  Indeed, in  some instances we may need to take one of these  graphs appearing partway, do another operation somewhere else, and then reverse the splits  done at $u$.  In order to keep track of the relevant information we will need to allow loops {\color \col{(which must have been created by splits at $u$)}} and special edges {\color \col{(to keep track of the edges created by splits at $u$)}}.  The following definition captures the basic properties these partway graphs will have.

\begin{definition}
Let $k\ge 4$ be an even number. A \emph{$k$-enhanced graph} is either an ordinary loopless graph or a graph $G$ with a unique \emph{special vertex} $u\in V(G)$ and a set of \emph{special edges}  $F\subseteq E(G)$ satisfying

\begin{itemize}
	\item $F \cap \delta(u) = \emptyset$, and
	\item if $e$ is a loop edge in $G$, then $e \in F$ (we may call $e$ a \emph{special loop}), and
	\item $d(u) + 2|F| = k$.
\end{itemize}	
\end{definition}

If $G$ is a $k$-enhanced graph with the special vertex $u$ and special edge set $F$, we associate with $G$ a  \emph{parent} graph $\overline{G}$ obtained from $G$ by subdividing each edge in $F$  and then identifying all of these newly created degree 2 vertices with the special vertex $u$.  Note that our assumptions imply that $\overline{G}$ will always be a loopless graph for which $u$ has degree $k$.

\bigskip

Observe that if $k\ge 4$ is even and $G$ is a graph that is obtained from a loopless graph $H$ by doing a number of splits at a vertex $u$ of $H$ with $d_H(u) = k$, then $G$ is a $k$-enhanced graph with $u$ as the special vertex and the set of edges created by the splits at $u$ as its set of special edges. The first two conditions follow from the assumption that $H$ is loopless and the last one is satisfied since $d_H(u) = k$. Hence, the definition above equips us with the appropriate terminology to talk about graphs that are partway through doing a complete split of a $k$-vertex $u$ of a loopless graph.

In the proof of our lemmas, there are times when we are dealing with such partway graphs and instead of doing the next split at $u$, we do another operation in $G$ to obtain $G'$ and then consider its parent graph  $\overline{G'}$. If the operation $O$ performed on $G$ to obtain $G'$ is such that $G'$ is also a $k$-enhanced graph, then $\overline{G'}$ will be a loopless graph in which $u$ has degree $k$.  In addition to this condition, if the graph $G$ is $k$-edge-connected except for possibly $\delta_{G} (u)$ and the operation $O$ is such that the $k$-enhanced graph $G'$ is also $k$-edge-connected except for possibly $\delta_{G'} (u)$, then observe that $\overline{G'}$ will be a $k$-edge-connected loopless graph. If these conditions are met and moreover, $G$ immerses a graph $K$ and the operation $O$ is such that $G'$ also immerses $K$, then $\overline{G'}$ is a $k$-edge-connected loopless graph which immerses $K$. The first condition motivates the definition of valid operations, adding the the second condition motivates the definition of safe operations, and finally adding the last one motivates the definition of good operations, introduced respectively in \cref{sub-valid}, \cref{sub-safe} and \cref{sub-good}.

\subsection{Valid Operations}
\label{sub-valid}

Let $G$ be a $k$-enhanced graph with the special vertex $u$ and special edge set $F$.  In order to maintain the essential features of this enhanced graph, we will insist upon using only certain operations.  The following definition captures these features.

\begin{definition}
If $G$ is a $k$-enhanced graph with no special vertices, then any split or deletion is \emph{valid}.  If $G$ has a special vertex $u$ and a special edge set $F$, then the following are defined to be \emph{valid} operations:
\begin{itemize}
\item deleting an edge not in $F \cup \delta(u)$;
\item splitting a 2-edge-path $xuy$ with $x, y \neq u$ to form a new edge $e$ and adding $e$ to $F$; 
\item splitting a 2-edge-path $xyz$ with $y \neq u$ and $| \{xy, yz\} \cap ( \delta(u) \cup F)| \le 1$.  If $F \cap \{xy,yz\} \neq \emptyset$, then we modify $F$ by replacing this edge by the newly created one.
\end{itemize}
Observe that the result of any of these operations is another $k$-enhanced graph with the special vertex $u$ and special edge set $F$ modified as indicated.

Accordingly, if $G$ has a special vertex $u$ with special edge set $F$, we call deleting an edge in $F \cup \delta(u)$ an \emph{invalid} edge deletion. Also, we call the split of a 2-edge-path $xyz$ with $y \neq u$ and $\{xy, yz\} \subseteq \delta(u) \cup F$ an \emph{invalid} split.

\bigskip

To help the reader remember the valid operations, we have illustrated them in \cref{validops}. In these illustrations, bold edges represent special edges, while non-bold edges represent non-special edges.

\begin{figure}[htbp]
	\centering
	\includegraphics{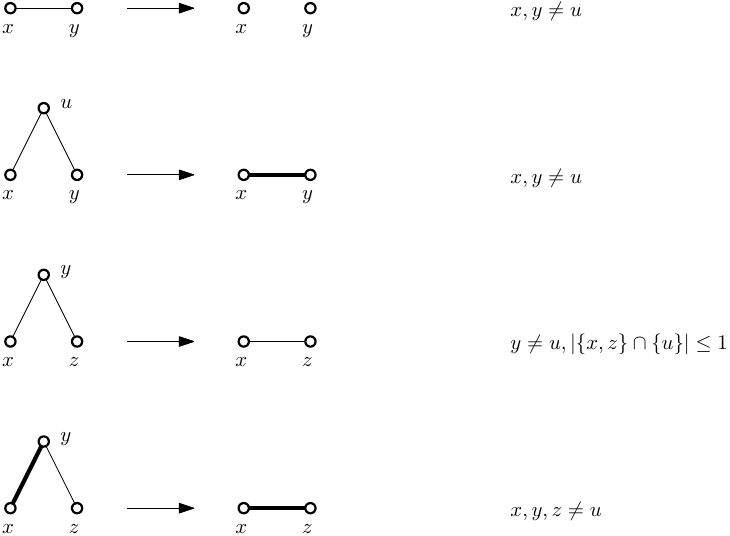}\\
	\caption{Valid operations}
	\label{validops}
\end{figure}

\end{definition}

\begin{note}
{\bf Throughout the rest of \cref{kecon}, we will assume that $\mathbf{k}$ is even and $\mathbf{ k\ge 4}$.}
\end{note}

\begin{notation}
	Let $O$ be any of the operations of edge deletion, splitting a 2-edge-path, or completely splitting a vertex of a graph (or $k$-enhanced graph) $G$. We denote by $O(G)$ the graph (or $k$-enhanced graph) obtained from $G$ by applying $O$.
\end{notation}

\subsection{Safe operations}
\label{sub-safe}

Next we will introduce the appropriate notion of edge-connectivity for $k$-enhanced graphs.  We say that a $k$-enhanced graph $G$ is \emph{nearly} $k$-\emph{edge-connected} if $G$ is a loopless $k$-edge-connected graph, or $G$ has a special vertex $u$ and special edge set $F$ and $G$ satisfies:
\begin{itemize}
\item $d(X) \ge k$ whenever $\emptyset \neq X \subset V(G)\setminus \{u\}$
\end{itemize}

Observe that whenever $G$ is a nearly $k$-edge-connected $k$-enhanced graph with a special vertex, the parent graph $\overline{G}$ will be $k$-edge-connected.  We are interested in performing valid operations that preserve our connectivity and this is captured by the following definition.

\begin{definition}
Let $G$ be a nearly $k$-edge-connected $k$-enhanced graph. We define a \emph{safe operation} on $G$ to be one of:
\begin{itemize}
\item A valid deletion of an edge
\item A valid split at the special vertex or a vertex of degree at least $k+2$,
\item A complete split of  either the special vertex or a vertex of degree $k$ where all splits are valid
\end{itemize}
with the additional property that the resulting $k$-enhanced graph is still nearly $k$-edge-connected.
Note that if $G$ has no special vertex, then the resulting graph will also have no special vertex (so must be another  $k$-edge-connected graph).  
\end{definition}

\cref{Mader} immediately gives us the existence of certain safe splits as follows.

\begin{lemma}
	\label{safesplit}
	Let $G$ be a nearly $k$-edge-connected $k$-enhanced graph with the special vertex $u$. Then the following hold.
	\begin{enumerate}[label=(\arabic*)]
		\item
		\label{safe-s}
		A safe split at the special vertex exists.
		\item
		\label{safe-ns}
		If $x$ is a vertex with $d(x) \ge k+2$ and $e(x,u) \le \frac{1}{2} d(u)$, then there exists a safe split at $x$.  
	\end{enumerate}
\end{lemma}

\begin{proof}
	For part \ref{safe-s}, note that thanks to the nearly $k$-edge-connectivity of $G$, $G$ is $d(u)$-edge-connected. Since $d(u)$ is even, it follows that $G$ does not have a cut-edge, so we can apply \cref{Mader} to do a split at $u$.  The resulting enhanced graph will be nearly $k$-edge-connected, as desired (Note that if the split creates a loop, it will be a special loop and thus we do not remove it.)   
	
	For part \ref{safe-ns}  let $F$ denote the set of special edges of $G$ and observe that $| \delta(x) \cap (F \cup \delta(u))| \le \frac{1}{2} d(u) + |F| \le \frac{1}{2}k$.  Therefore we may apply \cref{Mader} to split a 2-edge-path $yxz$ where $xy \not\in F \cup \delta(u)$.  The resulting enhanced graph will be nearly $k$-edge-connected, as desired.
\end{proof}

The next lemma concerns the existence of a safe complete split.

\begin{lemma}
	\label{safe-comp-sp}
	Let $G$ be a nearly $k$-edge-connected $k$-enhanced graph, and let $u'$ be a non-special vertex of degree $k$. Then a safe complete split of $u'$ exists.
\end{lemma}

\begin{proof}
	If $G$ does not have a special vertex, then we are done by \cref{Mader}, since it guarantees the existence of a complete split of $u'$ which preserves $k$-edge-connectivity. So a safe complete split of $u'$ exists. 
	
	Now, suppose $G$ has a special vertex $u$ and let $F$ denote the set of special edges of $G$.  Choose $O$ to be a complete split of $u'$, followed by removing any non-special loops, so that $O(G)$ is nearly $k$-edge-connected, and subject to this, the number of invalid splits by $O$ is minimum. Note that by Mader's Theorem, $O$ is well defined; let $G'= O(G)$ and let $F'$ denote the set of special edges in $G'$. Note that if no invalid split is done by $O$ then $G'$ is a nearly $k$-edge-connected $k$-enhanced graph, as desired.
	
	We claim that $O$ performs no invalid split. Otherwise, there are edges $xu', yu'\in F \cup \delta_G(u)$ that are split under $O$. 	 Since $G$ is nearly $k$-edge-connected and $d_G(u') = k$ we must have $e_G(u,u') \le \frac{d(u)}{2}$ (as otherwise $d_G( \{u,u'\}) < k$).  It follows from this that $| (F \cup \delta_G(u) ) \cap \delta_G(u')| \le \frac{k}{2}$.  Therefore the existence of $xu',yu' \in F \cup \delta_G(u)$ imply that there also must exist edges $zu', wu' \in \delta_G(u') \setminus (\delta_G(u)\cup F)$ so that $zu',  wu'$ are split under  $O$.  We claim that one of the following complete splits of $u'$,
	
	\begin{itemize}
		\item $O_1$ that splits the 2-edge-paths $xu'z$ and $yu'w$, and agrees with $O$ on other 2-edge-path splits,
		\item 
		$O_2$ that splits the 2-edge-paths $xu'w$ and $yu'z$, and agrees with $O$ on other 2-edge-path splits,
	\end{itemize} 
	contradicts the choice of $O$. Otherwise, since $O_1$ performs fewer invalid splits than $O$, the graph $G_1 = O_1(G)$ is not nearly $k$-edge-connected, and thus has an edge-cut $\delta_{G_1}(X)$, for some $X\subset V(G')$, other than $\delta_{G_1}(u)$ with size less than $k$. It follows from $d_{G'}(X)\ge k >d_{G_1}(X) $ together with $\delta_{G'}(X)- xy-zw \subseteq \delta_{G_1}(X)$ that $d_{G'}(X)= k$ or $k+1$, and that we may assume $\{x, z\}\subseteq X$ and $\{y, w\}\subseteq X^c$. A similar argument holds for $O_2$ and we may conclude that there exists $Y\subset V(G')$ with $d_{G'} (Y) \in \{k, k+1\}$ and that without loss of generality we may assume $\{x, w\}\subseteq Y$ and $\{y, z\} \subseteq Y^c$. In short, we have $x\in X\cap Y$, $z \in X\cap Y^c$, $w\in X^c \cap Y$, $y\in X^c \cap Y^c$.
	
	\bigskip
	\noindent \textbf{Remark.} Let $K$ be a graph with $Z \subset V(K), Z=Z_1 \cup Z_2, Z_1 \cap Z_2=\emptyset$. Then we have
	$$d(Z)=d(Z_1)+d(Z_2)-2e(Z_1, Z_2) \qquad  \qquad (*)$$

	Note that $d_{G'}(X) = d_{G'}(Y) = k+1$ is impossible. Otherwise, since $d_{G'}(X)$ is odd, by $(*)$ and possibly replacing $Y$ with $Y^c$, we may assume that $d_{G'}(X \cap Y)$ is even and $d_{G'}(X \cap Y^c)$ is odd. Similarly, since $d_{G'}(Y)$ is odd, we conclude that $d_{G'}(X^c \cap Y)$ is odd, and so both $X^c \cap Y$ and $X\cap Y^c$ contain a non-special vertex. So, nearly $k$-edge-connectivity of $G'$ implies that $d_{G'}(X^c\cap Y)\ge k+1$ and $d_{G'}(X\cap Y^c) \ge k+1$. Then, since $xy\in E_{G'}(X\cap Y, X^c\cap Y^c)$, we get the following contradiction in $G'$:
	$$2k+2 = d(X)+d(Y) = d(X^c\cap Y)+ d(X\cap Y^c) +2 e(X\cap Y, X^c\cap Y^c)\ge k+1+k+1+2.$$
	So, without loss of generality, we may assume $d_{G'}(X) =k$. Since none of $z, w$ are the special vertex, we have $d_{G'}(X^c\cap Y), d_{G'}(X \cap Y^c) \ge k$. Now, in $G'$ we get
	$$2k+1 \ge d(X)+ d(Y) = d(X^c\cap Y)+ d(X\cap Y^c) +2 e(X\cap Y, X^c\cap Y^c)\ge k+k+2,$$
	which is impossible. This contradiction completes the proof of the lemma.
\end{proof}

Next, we highlight situations where a safe edge deletion is possible:

\begin{lemma}
\label{lem1}
 Let $G$ be a nearly $k$-edge-connected $k$-enhanced graph. Suppose there exists $X \subset V(G)$ such that $ d(X) = k$, and every $x \in X$ has degree $k+1$.  Then there exists an edge in $G[X]$  
 whose deletion is safe.
\end{lemma}

\begin{proof}
Choose $Z\subseteq X$ such that $d(Z)=k$ and subject to this $Z$ is minimal. Since every $z \in Z$ has degree $k+1$, we have $|Z|\neq 1$, and $Z$ does not contain the special vertex,  if $G$ has one. If $G$ has a special vertex, call it $u$, and call the set of special edges $F$. Note that the degree assumption together with $|F| \le \frac{k}{2}$ imply that there is a non-special edge $e$ in $G[Z]$, so $e \in E(G[Z])\setminus (F \cup \delta(u))$.  To show that deleting $e$ is safe, it suffices to prove that $G\setminus e$ is nearly $k$-edge-connected. For a contradiction, suppose $e$ is in a $k$-edge cut $\delta (Y)$.
Note that $d(Z^c)=d(Z)=k$ implies that $Z^c$ contains at least one vertex $v$ other than the special vertex. We may assume (by possibly replacing $Y$ by $Y^c$) that $v\notin Y$. As $G$ is nearly $k$-edge-connected, $d(Z \cap Y), d({Z}^c \cap Y^c)\geq k$. However, it follows from
$$k+k\leq d(Z \cap Y) +d(Z^c \cap Y^c)\leq d(Z)+d(Y)=k+k$$
that $d(Z\cap Y)=k$, which contradicts minimality of $Z$.
\end{proof}

\begin{lemma}
\label{lem13}
Let $G$ be a nearly $k$-edge-connected $k$-enhanced graph, and suppose there exists $X \subset V(G)$ with $|X|\ge 2$ such that $ d(X)=k+1$, and every vertex in $X$ has degree $k+1$.  Then there exists an edge in $G[X]$ whose deletion is safe.
\end{lemma}

\begin{proof}
If $G$ has a special vertex, call it $u$, and let $F$ denote the set of special edges. Choose a non-special edge $e \in G[X]$, which has to exist because $|X|\geq 2$ and $|F| \le \frac{k}{2}$.  Note that $u \not\in  X$  (by the degree assumption) so $e  \not\in F \cup \delta(u)$.  So, if $G\setminus e$ is nearly $k$-edge-connected, then deletion of $e$ is safe. So, we may now assume that $e$ is in a $k$-edge-cut $\delta (Y)$.

Since $d(X)$ is odd, using $(*)$ (in the proof of \cref{safe-comp-sp}), by possibly replacing $Y$ with $Y^c$, we may assume that $d(X \cap Y^c)$ is even and $d(X \cap Y)$ is odd. Also, by $(*)$, we conclude that $d(X^c \cap Y)$ is odd, so $X^c \cap Y$ is nonempty and contains a non-special vertex (since the special vertex, if present, has even degree). Thus, both $X\cap Y^c$ and $X^c \cap Y$ contain a non-special vertex (note that one endpoint of $e$ is in $X^c \cap Y$). Now, by using \cref{cuts-edge-count} we have
$$2k+1=d(X)+d(Y)\geq d(X \cap Y^c) + d(X^c \cap Y) \geq 2k, $$ 
so by parity $d(X\cap Y^c)=k$. Therefore,  $\delta (X \cap Y^c)$ is a $k$-edge-cut where every vertex in $X\cap Y^c$ has degree $k+1$.  So, by applying \cref{lem1} we may conclude that there is an edge  in $G[X \cap Y^c]$ (and thus in $G[X]$) whose deletion is safe.
\end{proof}

\subsection{Good operations}
\label{sub-good}

Now we consider immersions for $k$-enhanced graphs. If $G$ is a $k$-enhanced graph  with a special vertex, we only say that $G$ immerses a graph $H$ if there exists such an immersion not using the special vertex as a terminal.  Building upon the previous subsection, we now embark on showing that the operations in the statement of \cref{sp-thm1} can be performed not only in a way that preserves near $k$-edge-connectivity, but also in a way that ensures that the resulting graph still contains an immersion of a particular $k$-edge-connected graph. This notion is formalized in the following definition.

\begin{definition}
	Let $G$ be a $k$-enhanced graph, and let $H$ be a $k$-edge-connected graph, with $G\succ H$. We define a  \emph{good operation} to be a {\color \col {safe operation (i.e. either}} a  safe split, a safe complete split, or a safe deletion of an edge) on $G$ which preserves an immersion of $H$ in the resulting $k$-enhanced graph.
\end{definition}

\begin{notation}
	If $G$ is a graph with $X \subset V(G)$, we denote the graph obtained by identifying $X$ to a single node by $G.X$.
\end{notation}

Let us make an observation before proceeding.

\begin{observation}
	Let $G$ be a graph with $X \subset V(G)$ such that there exists an immersion of $H$ with all terminals in $X$. Then $G.X^c$ contains $H$ as an immersion.
\end{observation}

We begin with one lemma where a good split is easily found.

\begin{lemma}
\label{overloadedspecial}
Let $G$ be a nearly $k$-edge-connected $k$-enhanced graph with a special vertex $u$ and another vertex $x$ satisfying $e(u,x) > \frac{1}{2}d(u)$. Then there exists a good complete split of $u$.
\end{lemma}

\begin{proof}
Choose a complete split of $u$ so that every edge incident with $u$ but not $x$ is split off with an edge between $u$ and $x$. Note that since $e(u,x) > \frac{1}{2}d(u)$, some special loops will be formed at $x$. Delete any loops formed in this process and let $G'$ be the resulting graph.  Since $G'$ is isomorphic to the graph obtained from $G$ by identifying $u$ and $x$, it is $k$-edge-connected. Also, since $G$ has an immersion of $H$ not using the vertex $u$ as a terminal, $G'$ must still immerse $H$, so we have found a good complete split of $u$.
\end{proof}

The following lemma enables us to handle $k$-edge-cuts:

\begin{lemma}
\label{lem2}
Let $G$ be a nearly $k$-edge-connected $k$-enhanced graph. If $G$ has a nontrivial $k$-edge-cut $\delta (X)$ such that some immersion of $H$ has no terminal in $X$, then there exists a good operation.
\end{lemma}

\begin{proof}
If $G$ has a special vertex, we will call it $u$, and let $F$ denote the set of special edges. If the special vertex $u$ is in $X$, we will do a safe split at $u$ using \cref{safesplit}\ref{safe-s}.  Else, thanks to \cref{overloadedspecial}, we may assume for every vertex $x \neq u$ we have $e(u,x) \le \frac{1}{2} d(u)$. So, if there exists a vertex $x \in X$ with $d(x)\ge k+2$, we can use \cref{safesplit}\ref{safe-ns} to do a safe split at $x$. Also, if there exists a vertex $x \in X$ with $d(x)=k$, we use \cref{safe-comp-sp} to do a safe complete split at $x$. In the only remaining case, every vertex in $X$ is of degree $k+1$, and by applying \cref{lem1}, we can do a safe deletion of an edge in $G[X]$.

In either case, let $G'$ be the resulting $k$-enhanced graph. To see that the operation is good, observe that there remain $k$ edge-disjoint paths in $G'$ between any pair of non-special vertices, one in $X$, and the other in $X^c$ (since $G'$ is (nearly) $k$-edge-connected). Therefore $G'$ immerses $G.X$, and thus it immerses $H$.
\end{proof}

The next two lemmas which concern a broader family of graphs, will later be helpful in dealing with $(k+1)$-edge-cuts in $k$-enhanced graphs.

\begin{lemma}
\label{lem3}
Let $G$ be an internally $k$-edge-connected graph in which every vertex of degree less than $k$ has even degree. If $d(x)$ is odd, then there exists $y \in V(G)\setminus x$ such that $\lambda (x, y)\geq k+1$.
\end{lemma}

\begin{proof}
We prove the statement by induction on $|V(G)|$. Note that by parity, there must exist another vertex $y \in V(G)$ with odd degree. Note that since $d(y)$ is odd and $k$ is even, we have $d(y)\ge k+1$. If every edge-cut separating $x$ from $y$ has size at least $k+1$, by Menger's Theorem (\cref{Menger}) we are done. Otherwise, there exists an edge-cut $\delta (Y)$ {\color \col of size at most $k$}, with $y \in Y$, separating $x$ from $y$. Note that since $d(y)\ge k+1$, we have $|Y|\geq 2$ and so $d(Y) = k$.

{\color \col Since} $|Y|\geq 2$, the graph ${G'}=G.Y$ which satisfies the lemma's hypothesis, has fewer vertices than $G$. Also $x$ has odd degree in $G'$ as well, thus by induction hypothesis there exists $y' \in V(G')\setminus x$ such that $\lambda _{G'}(x, y')\geq k+1$. It follows, however, that $\lambda _{G}(x, y')\geq k+1$ as well, since $\lambda_G(x,y)=k$ implies that $G \succ  {G'}$.
\end{proof}

\begin{lemma}
\label{lem4}
Let $G$ be an internally $k$-edge-connected graph in which every vertex of degree less than $k$ has even degree.
If $\delta (X)$ is a $(k+1)$-edge-cut in $G$, there exist $x \in X, y \in X^c$ such that $\lambda(x,y)\geq k+1$.
\end{lemma}

\begin{proof}
Let $G_1 = G.X$ and $G_2 = G.X^c$, with $s$ and $t$ being the nodes replacing $X$ and $X^c$, respectively. Note that both $G_1$ and $G_2$ satisfy the hypothesis of \cref{lem3}. Since $s$ is a vertex of odd degree in $G_1$, by \cref{lem3}, there exists $y \in X^c$ such that $\lambda_{G_1}(s, y)\geq k+1$. Note that this also implies that by taking $y$ to be the terminal corresponding to $t$, we have $G \succ G_2$. It can be similarly argued that there exists $x\in X$ such that $\lambda_{G_2}(x,t)\geq k+1$ and since by taking $y$ to be the terminal corresponding to $t$, we get an immersion of $G_2$ in $G$, we have $\lambda_G(x,y)\geq k+1$.
\end{proof}

Having the lemma above in hand, we can now efficiently handle $(k+1)$-edge-cuts in $k$-enhanced graphs:

\begin{lemma}
\label{lem5}
 Let $G$ be a nearly $k$-edge-connected $k$-enhanced graph. Suppose $G$ has a nontrivial $(k+1)$-edge-cut $\delta (X)$ such that some immersion of $H$ has no terminal in $X$, then there exists a good operation.
\end{lemma}

\begin{proof}
If $G$ has a special vertex, we will call it $u$, and let $F$ denote the set of special edges. If $u$ is in $X$, we will do a safe split at $u$ using \cref{safesplit}\ref{safe-s}.  Else,  by \cref{overloadedspecial}, we may assume for every non-special vertex $x$ we have $e(u,x) \le \frac{1}{2} d(u)$. So, if there exists a vertex in $X$ of degree other than $k+1$, we will apply \cref{safesplit}\ref{safe-ns} or \cref{safe-comp-sp} to either do a safe split or do a safe complete split at such a vertex. In the remaining case, every vertex in $X$ has degree $k+1$, and 
we can apply \cref{lem13} to do a safe edge deletion from $G[X]$.

In any case, we claim the safe operation is good. Let $G'$ be the $k$-enhanced graph resulting from doing the safe operation. Note that $\delta (X)$ remains a $(k+1)$-edge-cut in $G'$, since doing a split changes the size of an edge-cut by an even number, and by the (nearly) $k$-edge-connectivity of $G'$, $d_{G'}(X) \geq k$. We may now apply \cref{lem4} to choose $x \in X$, $y \in X^c$ with $\lambda(x,y) \geq k+1$. Thus $G'$ immerses $G.X$, and therefore it immerses $H$.
\end{proof}

\subsection{Proof of the main theorem}
\label{sec-pf}

In this subsection, we use the results in the last two subsections to prove \cref{sp-thm1}.  We begin by recording a basic observation that follows immediately from our definitions.

\begin{observation}
\label{sew}
Let $G$ and $H$ be $k$-edge-connected graphs with $G \succ H$. Suppose $u \in V(G)$ has degree $k$. Let $K$ be a nearly $k$-edge-connected $k$-enhanced graph with the special vertex $u$ obtained from $G$ by doing splits at $u$ (note that any such split creates a special edge) and the special edge set $F$. Assume that $K \succ H$ and that there exists a good operation apart from a single split at $u$ which takes $K$ to $K'$.  Then there is a corresponding good operation taking $G$ to $G'$ so that $\overline{K'} \cong G'$. 
\end{observation}

\begin{proof}
	Note that by definition, $\overline{K} \cong G$. Let $O$ be the good operation taking $K$ to $K'$. If $O$ interacts only on edges other than $F \cup \delta_K(u)$, i.e. if $O$ is either a safe edge deletion, or a safe (complete) split at a vertex other than $u$ which does not involve any edge in $F$, then by letting $G' = O(G)$, we have $\overline{K'} = \overline{O(K)}  \cong O(\overline{K}) \cong O(G) = G'$.
	
	Now, suppose $O$ is a safe split of a 2-edge path $xyz$ that involves special edges. Since $O$ is a valid split, without loss of generality we may assume that $xy \in F$ and $yz \notin F \cup \delta_K(u)$. Now, if we let $O'$ be splitting the 2-edge-path $uyz$ in $G$ and let $G' = O'(G)$, then observe that $\overline{K'} \cong G'$. It is straightforward to see that $O'$ is a good operation.
\end{proof}

The next three lemmas concern the three operations allowed in stepping from $G$ towards $H$, which are $k$-edge-connected graphs with $G\succ  H$, and show that in each case we can take a step maintaining $k$-edge-connectivity.

\begin{lemma}
\label{lem6}
Let $G$ and $H$ be loopless $k$-edge-connected graphs with $G\succ  H$. If there is a complete split of a $k$-vertex $u$ of $G$ preserving an $H$-immersion, then there is a good operation.
\end{lemma}

\begin{proof}
Consider the complete split of $u$ as $\frac{k}{2}$ many splits at $u$, and choose a sequence of splits which, while preserving an $H$-immersion, results in the fewest number of loops. If $u$ could be completely split without ever creating a too small of an edge-cut other than $\delta(u)$ along the way, we are done. Otherwise, we will stop doing these splits the first time the resulting graph $K$ is about to
have an edge-cut of size less than $k$ other than $\delta(u)$. In $K$, therefore, there exists a subset $X \neq \{u\}, \{u\}^c$ of $V(G)$ for which $d_{K}(X)\geq k $; however, doing the next split makes the size of its boundary to drop to less than $k$, so $d_{{K}}(X)=k$ or $k+1$. Now, let $F$ be the set of edges, including loops,  of $K$ which are created so far by doing splits at $u$. Note that since $G$ is loopless and every edge in $F$ is created by a single split at $u$, we have $F\cap \delta_{K}(u) = \emptyset$, every loop is in $F$, and $2|F| + d_{K}(u) = k$, so $K$ together with $u$ as the special vertex and $F$ as the set of special edges is a  nearly $k$-edge-connected  $k$-enhanced graph.

Now, since completely splitting $u$ results in $d(X)<k$ and preserves an immersion of $H$, there is an immersion of $H$ with all terminals on one side of $\delta (X)$, say $X^c$. First, suppose $\delta (X)$ is a nontrivial cut.
Then \cref{lem2} or \cref{lem5} applied to $K$ guarantee the existence of a good operation
 $O$, which may or may not be a (complete) split at $u$. Recall that by definition of a good operation, $O(K)$ is a nearly $k$-edge-connected  $k$-enhanced graph which immerses $H$. Now, if $O$ is a (complete) split at $u$, we resume splitting $u$ in $O(K)$ and else, we let $G^*= \overline{O(K)}$.  By \cref{sew}, we know that $G^*$ may be obtained from $G$ by a good operation, so by definition of a good operation $G^*$ is $k$-edge-connected and $G \succ H$, as desired.  

Now suppose $\delta(X)$ is a trivial cut with $X=\{ v \}$. Therefore the next split at $u$ would create a loop at $v$. Note that there cannot be a vertex $w \in N_{K}(u)\setminus v$, because if there was one, then we could have split off $vuw$ to get the graph $K'$ instead. It is because $K'$ has no loop and immerses $H$, since splitting off $wvu$ in $K'$ results in the same graph as splitting off $vuv$ in $K$ would.

Therefore $N_{K}(u)=\{v\}$, implying that in $K$, $d(v)=d(\{u, v\}) +d(u)$. This, however, contradicts $d(X)=k$ or $k+1$, as $d(X)=d(v)=d(\{u, v\}) +d(u)\ge k+ d(u)\ge k+2$, where the inequalities hold because $K$ is nearly $k$-edge-connected, and $u$ is of even degree. This completes the proof.
\end{proof}

\begin{lemma}
\label{lem7}
Let $G$ and $H$ be $k$-edge-connected graphs with $G \succ H$. If there is an edge $e$ such that $G \setminus e$ has an $H$-immersion, then a good operation exists.
\end{lemma}

\begin{proof}
Suppose $e$ is in a $k$-edge-cut. If it is incident with a $k$-vertex $u$, then by previous lemma, there is a good operation. Otherwise, $e$ is in a nontrivial $k$-edge-cut, with all terminals of $H$ on one side of the cut, thus we can use \cref{lem2} to find a good operation.
\end{proof}

\begin{lemma}
\label{lem8}
Let $G$ and $H$ be $k$-edge-connected graphs with $G \succ H$. If there is a split at a vertex $v$ preserving an $H$-immersion, then a good operation exists.
\end{lemma}

\begin{proof}
Suppose splitting at $v$ makes an edge-cut $\delta (X)$ too small, then $d(X)=k$ or $k+1$. Also, all terminals of $H$ are on one side of the cut, say $X^c$. If $\delta (X)$ is a nontrivial edge-cut \cref{lem2} or \cref{lem5} may be applied. If $|X|=1$, with $d(X)=k$, we apply \cref{lem6} to completely split the vertex in $X$, and if $d(X)=k+1$ we will apply \cref{lem7} to delete an edge incident to it.
\end{proof}

The proof of \cref{sp-thm1} is now immediate:

\begin{proof}[Proof of \cref{sp-thm1}.]
		Since $H \prec G$, there is either a complete split at a vertex of degree $k$, or an edge deletion or a split at a vertex of degree at least $k+2$ that takes $G$ to $G'$ such that $G'\succeq H$. Now, apply  \cref{lem6}, \cref{lem7}, or \cref{lem8}.
\end{proof}
\section{3-edge-connected, internally 4-edge-connected graphs}
\label{in4econ}

In this section we establish \cref{sp-thm2intro}. Later, as an application, we will see that if a $3$-edge-connected, internally $4$-edge-connected graph other than $K_{2,2,2}$ immerses $K_5$, it also has a $K_{3,3}$-immersion. First, we move towards proving \cref{sp-thm2intro}, which, for convenience is {\color \col{repeated}} here.
\begin{theorem}
\label{sp-thm2}
Let $G$ and $H$ be $3$-edge-connected and internally $4$-edge-connected loopless graphs, with $G \succ H$. Assume that $|V(H)|\ge 2$, and $(G, H) \ncong (Q_3, K_4), (Q_3, K_2^3)$. Then there exists an operation taking $G$ to $G'$ such that $G'$ is $3$-edge-connected, internally $4$-edge-connected, and $G' \succeq   H$, where an operation is either
\begin{itemize}
\item deleting an edge,
\item splitting at a vertex of degree at least 4,
\end{itemize}
each followed by iteratively deleting any loops, and suppressing vertices of degree $2$.
\end{theorem}
As in the proof of \cref{sp-thm1}, we will consider each operation separately, and the proof of the theorem will then be immediate. First, we will modify our notion of a good operation from \cref{kecon} as follows:

\begin{definition}
Let $G$ and $H$ be 3-edge-connected and internally 4-edge-connected loopless graphs, with $G \succ H$. We define a \emph{good operation} to be either a split at a vertex of degree at least 4, or a deletion of an edge from $G$ which preserves 3-edge-connectivity, internal 4-edge-connectivity, and an immersion of $H$ in the resulting graph.
\end{definition}
\begin{lemma}
\label{lem9}
Let $G$ and $H$ be as in \cref{sp-thm2}. Suppose there is an edge $e \in E(G)$ such that $G \setminus e$ has an $H$-immersion. If $(G, H) \ncong (Q_3, K_4), (Q_3, K_2^3)$, then a good operation exists.
\end{lemma}
\begin{proof}
Since deletion of $e$ is followed by suppression of any resulting vertices of degree two, and as $G$ is internally 4-edge-connected, $G\setminus e$ is clearly $3$-edge-connected. If deletion of $e$ does not preserve internal 4-edge-connectivity, then $e$ must be contributing to some $4$-edge-cut $\delta (X)$ in which each side has either at least three vertices, or has two vertices which are not both of degree 3. We call such a cut an \emph{interesting cut}.

Note that $H$ too is internally 4-edge-connected, thus all but possibly one, of the terminals of an immersion of $H$ lie on one side of this cut, say $X$. Let $X'$ be the maximal subset of $V(G)$ containing $X$, such that $\delta (X')$ is interesting. Suppose there is an edge $uv$ in ${X'}^c$ not contributing to an interesting edge-cut, then deleting $uv$ is a good operation. It is because $G \setminus uv$ is $3$-edge-connected, internally $4$-edge-connected. Also $G \setminus uv$ has an $H$-immersion, because it immerses $(G\setminus e).{X}^c$.

We may now assume that $uv$ is in some interesting edge-cut $\delta (Y)$. Note that maximality of $X'$ implies that $X' \cap Y, X' \cap Y^c \neq \emptyset $. Because if, say, $X' \cap Y = \emptyset$ then $X' \subsetneq Y^c$ and since $\delta(Y)$ is an interesting edge-cut, we get a contradiction with maximality of $X'$. Also, we claim that there cannot be edges contributing to both $\delta (X'), \delta(Y)$. To prove the claim, suppose, to the contrary, that there are edges between, say, $X' \cap Y$ and ${X'}^c \cap Y^c$, i.e. $e(X' \cap Y, {X'}^c \cap Y^c) \neq 0$; we will denote this quantity by $f$ (as in \cref{two4cutcross}).
\begin{figure}[htbp]
\centering
  \includegraphics[height=5cm]{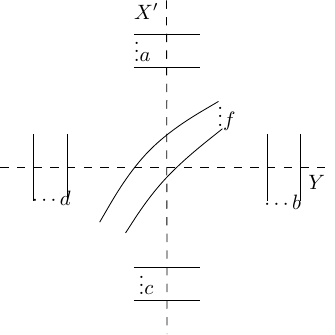}\\
  \caption{ Cuts $\delta (X')$, $\delta (Y)$ relative to each other}
  \label{two4cutcross}
\end{figure}
 Then it follows from
$$ 8=d(X')+d(Y) = d({X'}^c \cap Y) +d(X' \cap Y^c)+2f  \geq 3+3+2f$$
that if $f \neq 0$, it equals to 1 and moreover, $d({X'}^c \cap Y)=d(X' \cap Y^c)=3$.

Using a similar argument, one can see that, if in addition to $f \neq 0$, there were also edges between $X'\cap Y^c, {X'}^c\cap Y$,
 then $d({X'}^c \cap Y^c)=3$. Thus, both ${X'}^c \cap Y^c, {X'}^c \cap Y$ would consist of a single vertex of degree $3$, contradicting $\delta (X')$ being interesting. Therefore the number of edges contributing to both $\delta(X'), \delta (Y)$ equals $f$.

We will now show that $f \neq 0$ results in a contradiction. We define (please see \cref{two4cutcross})
\[a = e(X'\cap Y^c, {X'}^c \cap Y^c) \hspace{1cm} b = e( {X'}^c \cap Y^c,  {X'}^c \cap Y)\]
\[c = e( {X'}^c \cap Y, X' \cap Y) \hspace{1cm} d= e(X'\cap Y, X'\cap Y^c)\]
 Note that from $d({X'}^c \cap Y)=3$ we may conclude, without loss of generality, that $b \geq 2$. Now, by alternatively looking at the cuts $\delta ({X'}^c\cap Y)$, $\delta (X')$, $\delta (X'\cap Y^c)$, $\delta (Y) $, we see that if $b\ge 2$, then $c\le 1$, so $a \ge 2$, thus $d \le 1$. Therefore, $d(X'\cap Y)=c+ d+ f \le 3$, so $X'\cap Y$ consists of a single vertex of degree three. However, this, together with the earlier conclusion that $X'\cap Y^c$ consists of a single vertex of degree three, contradicts $\delta(X')$ being interesting. Therefore $f=0$, so there are no edges contributing to both $\delta(X')$ and $\delta(Y)$.

Now, we show that $a=b=c=d=2$. For a contradiction, we will assume that, say $a>2$, and, similar to the argument above, alternatively look at the cuts $ \delta (X'), \delta (X'\cap Y), \delta (Y) $. It then follows that $c\le 1$, so $d \ge 2$, thus $b\le 2$. So, in order for $d({X'}^c\cap Y)=b+c\ge 3$, we must have $b=2, c=1$. Also, we have $d(Y)=4=b+d$, so $d=2$, thus $d(X'\cap Y)=c+d=3$. 
Hence, each ${X'}^c\cap Y$ and $X'\cap Y$ consist of a single vertex of degree three, which contradicts $\delta(Y)$ being interesting.

Therefore, $a=b=c=d=2$, and thus $\delta({X'}^{c} \cap Y), \delta({X'}^c \cap Y^c)$ are 4-edge-cuts. However, by maximality of $X'$, they cannot be interesting cuts. Thus each of ${X'}^{c} \cap Y, {X'}^c \cap Y^c$ consists of only one vertex, or two vertices both of degree 3.

We are now ready to prove that a good operation exists unless $(G, H) \cong (Q_3, K_4)$ or $(G, H) \cong (Q_3, K_2^3)$. Consider different possibilities for ${X'}^{c} \cap Y, {X'}^c \cap Y^c$:

\begin{itemize}
\item Both sets consist of one vertex, see \cref{onevx}$(a)$. Here, a good operation is to split off $wuv$. Note that the resulting graph immerses $H$, as it immerses $(G\setminus e).{X}^c$.
\item Only one set consists of one vertex. Then it is easy to verify that ${X'}^c$ should be as in \cref{onevx}$(b)$. Here, deleting $vw$ is a good operation. 
\begin{figure}[htbp]
\centering
  \includegraphics[height=4.5cm]{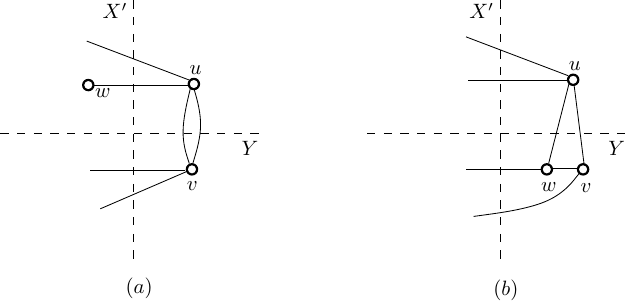}\\
  \caption{At least one of ${X'}^{c} \cap Y, {X'}^c \cap Y^c$ consists of only one vertex}
  \label{onevx}
\end{figure}

\item Both sets have two vertices in them, see \cref{del31}.
\begin{figure}[htbp]
\centering
  \includegraphics[height=4cm]{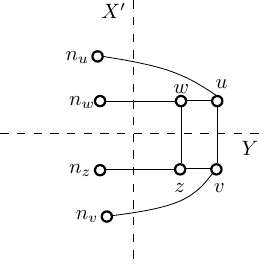}\\
  \caption{Both ${X'}^{c} \cap Y, {X'}^c \cap Y^c$ consist of two vertices}
  \label{del31}
\end{figure}
Here the operation will be deleting $uw$ or $vz$, from which we claim at least one is a good operation unless $G \cong Q_3$. Suppose that deleting both $uw$ and $vz$ destroy internal 4-edge-connectivity, thus both these edges contribute to some interesting cuts.

As before, it can be argued that the cuts look like as in  \cref{del32} with respect to each other.
\begin{figure}[htbp]
\centering
  \includegraphics[height=4cm]{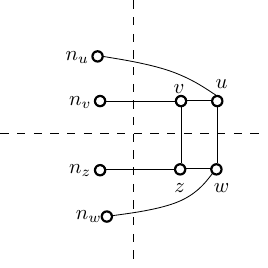}\\
  \caption{Both $uw$ and $vz$ are in interesting edge-cuts}
  \label{del32}
\end{figure}
Now, ignoring $\{ u, v, w, z\}$ in \cref{del31}, and \cref{del32}, we can see that there exists a 2-edge cut separating $\{ n_u, n_w\}$ from $\{n_v, n_z\}$, and another one separating $\{n_u, n_v\}$ from $\{ n_w, n_z\}$, implying that $n_u, n_v, n_w, n_z$ form a square, thus $G\cong Q_3$. It has now only remained to notice that $K_4, K_2^3$ are the only internally 4-edge-connected graph that $Q_3$ immerses.
\end{itemize}
\end{proof}

Our next task is to deal with splits in $G$ that preserve an $H$-immersion, which will be done in \cref{lem10}. The following statement, which holds for a broader family of graphs than the ones we work with, features in the proof of \cref{lem10}.
\begin{lemma}
\label{3econdel}
Let $H$ be a $3$-edge-connected graph, and let $Y$ be a minimal subset of $V(H)$ such that $\delta(Y)$ is a nontrivial $3$-edge-cut in $H$. Then for every edge $e$ in $H[Y]$, $H\setminus e$ is internally $3$-edge-connected.
\end{lemma}
\begin{proof}
For a contradiction, suppose an edge $e=yz$ in $H[Y]$ contributes to some nontrivial $3$-edge-cut $\delta(Z)$, where $z\in Z$. We will look into how $Y,Z$ look like with respect to one another. Note both $Y\cap Z$ and $Y\cap Z^c$ are nonempty, as $z\in Y\cap Z, y\in Y\cap Z^c$. Also, both $Y^c\cap Z$ and $Y^c\cap Z^c$ are nonempty. It is because, if, say $Y^c\cap Z=\emptyset$, then  $\delta(Y\cap Z)$ would be a nontrivial $3$-edge-cut, which contradicts the choice of $Y$, as $Y\cap Z\subsetneq Y$.

Now, since $H$ is $3$-edge-connected, we have $d(Y\cap Z),d(Y^c\cap Z^c)\ge 3$. It now follows from
$$6 \le d(Y\cap Z)+d(Y^c\cap Z^c)+2e(Y^c\cap Z, Y\cap Z^c)=d(Y)+d(Z) = 6$$
 that $d(Y\cap Z)=d(Y^c\cap Z^c)=3$ and $e(Y^c\cap Z, Y\cap Z^c)=0$. Similarly, we obtain $d(Y\cap Z^c)=d(Y^c\cap Z)=3$ and $e(Y\cap Z, Y^c\cap Z^c)=0$. Now, since
 $$d(Y)=3=e(Y\cap Z, Y^c\cap Z)+e(Y\cap Z^c,Y^c\cap Z^c),$$
  we have, say, $e(Y\cap Z, Y^c\cap Z)\le 1$. Similarly, it follows from $d(Z)=3$ that we have, say, $e(Y\cap Z, Y\cap Z^c)\le 1$. Hence, $d(Y\cap Z)\le 2$, a contradiction.
\end{proof}
\begin{lemma}
\label{lem10}
Let $G$ and $H$ be as in \cref{sp-thm2}. Assume that there is a split at a vertex $v \in V(G)$ preserving an $H$-immersion. If $(G, H) \ncong (Q_3, K_4), (Q_3, K_2^3)$, then a good operation exists.
\end{lemma}

\begin{proof}
Let $uvw$ be the $2$-edge-path that is to be split. Note if $d(v)=3$, then deleting the edge incident to $v$ other than $vu,vw$ preserves the $H$-immersion. Hence, by \cref{lem9} we are done. Also, observe that if a split is done at a vertex of degree at least four, the resulting graph is $3$-edge-connected. Therefore, we only need to look into the case where splitting off $uvw$ destroys internal $4$-edge-connectivity. So, it must be the case that $uv, vw$ contribute to some 4- or 5-edge-cut $\delta (X)=\{ uv, wv, x_1y_1, x_2y_2 (, x_3y_3): u, w, x_i \in X \}$, where $|X|,|X^c|\ge 2$. We now split the analysis into cases depending on $d(X)$.
\begin{claim}
\label{splitdX4}
If $d(X)=4$, a good operation exists.
\end{claim}
Since $H$ is $3$-edge-connected and splitting $uvw$ creates a 2-edge-cut, all terminals of $H$ lie on one side of the cut. Also, since $G$ is $3$-edge-connected, each side of the cut contains an edge lying in it, i.e. $E(G[X]), E(G[X^c])\neq \emptyset$.

First, suppose all terminals of $H$ are in $X$. Observe that if we can delete an edge in $G[X^c]$ in a way that  the connectivity between $y_1$ and $y_2$ in $G[X^c]$ is preserved, an $H$-immersion is present in the resulting graph. We propose to delete an edge $e \in E(G[X^c])$ and claim that deleting $e$ preserves the $H$-immersion. It suffices to show $e$ is not a cut-edge in $G[X^c]$ separating $y_1, y_2$. For a contradiction, suppose $e=\delta(Y)$ separates $y_1,y_2$ in $G[X^c]$, where $y_1\in Y$. We may also assume, without loss of generality, that $v \in Y$. Then $\delta_G(Y^c)$ would be a $2$-edge-cut in $G$, a contradiction. Therefore, we can delete $e$ using \cref{lem9}.

Next, suppose all terminals of $H$ are in $X^c$. Similar to the previous case, if we modify $X$ in a way that preserves the connectivity between $x_1, x_2$ in $G[X]$, then an $H$-immersion is sure to exist in the resulting graph. Again, we propose to delete an edge $e \in E(G[X])$ and claim that deleting $e$ preserves the $H$-immersion. It suffices to show $e$ is not a cut-edge in $G[X]$ separating $x_1, x_2$. For a contradiction, suppose $e=\delta(Y)$ separates $x_1,x_2$ in $G[X]$, where $x_1\in Y$. Note $3$-edge-connectivity of $G$ implies that $\delta(Y)$ separates $u, w$ as well. We may assume, without loss of generality, that $u\in Y, w\in Y^c$. Then $d_G(Y)=d_G(Y^c)=3$, thus it follows from internal $4$-edge-connectivity of $G$ that $|Y|=|Y^c|=1$ and $Y=\{u=x_1\}, Y^c=\{w=x_2\}$. Therefore, $X$ consists of two vertices $u,w$ of degree three, and thus deleting $uw$ preserves the $H$-immersion.

\begin{claim}
If $d(X)=5$, a good operation exists.
\end{claim}
By the internal edge-connectivity of $H$, all terminals of $H$ but possibly one, lie on one side of the cut. First, suppose that most terminals of $H$ are in $X$. Observe that if $X^c$ is modified in a way that preserves the presence of three edge-disjoint paths form a vertex in it to $X$ not using $uv, vw$, the presence of $H$-immersion is guaranteed. Next, suppose that most terminals of $H$ are in $X^c$. In this case, if we manage to modify $X$ in a way that preserves the presence of three edge-disjoint paths form a vertex in it to $X^c$ covering $\delta(X)\setminus \{uv,vw\}$, the presence of $H$-immersion is guaranteed. 

Now, we claim the aforementioned modifications are possible whether most terminals are in $X$ or in $X^c$. In any case, let $G'$ be the graph resulting from splitting off $uvw$, followed by suppressing $v$ in case $d_{G}(v)=4$. We denote the edge created by splitting $uvw$ by $e'$. Note, by \ref{splitdX4}, we may assume $G'$ is 3-edge-connected.

Take an arbitrary nontrivial $3$-edge-cut $\delta_{G'}(Y)$ in $G'$. Observe that $\delta_G(Y)$ must have been a $5$-edge-cut in $G$, which both edges of the split $2$-path $uvw$ contributed to. So, in particular, $e'$ lies either completely in $Y$ or in $Y^c$. Also, there must be an edge other than $e'$ in $G'[Y]$. It is because $3$-edge-connectivity of $G$ implies $6\le \sum_{y \in Y} d_G(y)=d_G(Y)+2 e_G(G[Y])=5+2 e_G(G[Y])$. Thus $e_G(G[Y])>0$, and so there is an edge $\neq e'$ in $G'[Y]$.

Now, let $Z$ denote the side of $\delta(X)$ containing most terminals of $H$ (, so $Z=X$ or $X^c$). We will show that there is an edge lying in $Z^c$ which we could delete, while preserving an $H$- immersion. Since $\delta_{G'}(Z)$ is a nontrivial $3$-edge-cut, we may choose a minimal $Y\subseteq Z^c$ such that $\delta_{G'}(Y)$ is a nontrivial $3$-edge-cut.

It is argued above that there exists an edge $e\neq e'$ in $G'[Y]$. We claim deletion of $e$ preserves the $H$-immersion. It is because it follows from \cref{3econdel} that $G'\setminus e$ is internally $3$-edge-connected. Now, $3$-edge-connectivity of $G'$ and $d_{G'}(Y)=3$ imply that $G'[Y]\setminus e$ has a vertex of degree at least three. Therefore, there exists in $G'\setminus e$ three edge-disjoint paths from such a vertex to $Z$. Observe that since these set of paths cover $\delta(Z)$, deletion of $e$ from $G$ preserves the presence of $H$-immersion. We now can use \cref{lem9} to delete $e$ from $G$.
\end{proof}

Proof of \cref{sp-thm2} is now immediate.

\begin{proof}[Proof of \cref{sp-thm2}]
	Since $H \prec G$, there is either a split  or an edge deletion that takes $G$ to $G'$ such that $H\preceq G'$. Now, apply  \cref{lem9} or \cref{lem10}.
\end{proof}

Having established \cref{sp-thm2}, we will now take advantage of it to establish a chain theorem for the family of $3$-edge-connected, and internally $4$-edge-connected graphs.

\begin{cor}
	\label{sp-chain-3ec-int4ec}
	Let $G$ be a $3$-edge-connected and internally $4$-edge-connected graph, where $|V(G)| \ge 2$. If $G\ncong Q_3, K_2^3$, one of the operations in the statement of \cref{sp-thm2} may be applied to $G$, where the resulting graph is $3$-edge-connected, internally $4$-edge-connected.
\end{cor}

\begin{proof}
	Since $G$ is $3$-edge-connected, it immerses $K_2^3$. Now, apply \cref{sp-thm2} for $H = K_2^3$.
\end{proof}

We now work towards proving \cref{corintro}. The idea is to examine $3$-edge-connected, internally $4$-edge-connected graphs ``one step bigger", or perhaps ``a few steps bigger", than $K_5$, and see if they immerse $K_{3,3}$. One subtlety here is that we are working with multigraphs, thus even graphs ``much bigger than" $K_5$ may happen to be on five vertices, and thus not possess $K_{3,3}$-immersions. Therefore, we need some tool to limit the graphs necessary to examine. Given that $K_5$ itself is $4$-edge-connected, \cref{lem11} serves very well in doing so. First, however, we need the following definition.

\begin{definition}
We define a \emph{good sequence} from $G$ to $H$ to be a sequence of graphs
$$G=G_l, G_{l-1}, \ldots , G_2, G_1, G_0\cong H$$
in which each $G_i$ is 3-edge-connected, and internally 4-edge-connected, and $G_i$ is resulting from applying an operation $o_{i+1}$ (as defined in the statement of \cref{sp-thm2}) to $G_{i+1}$.
\end{definition}

\begin{lemma}
\label{lem11}
Let $G$ be $3$-edge-connected, internally $4$-edge-connected, and $H$ be $4$-edge-connected. Suppose there is a good sequence from $G$ to $H$, and choose a good sequence from $G$ to $H$
$$G=G_l, G_{l-1}, \ldots , G_2, G_1, G_0\cong H$$
such that $\min \{ k: |V(G_k)|>|V(H)| \}$ is as small as possible. Then either
\begin{itemize}
\item[(a)] $G_1$ is as in \cref{gi}(a), with $v_1 \neq v_2$, $v_3 \neq v_4$, and the last operation, $o_1$, is to split off $v_1uv_2$ and $v_3uv_4$.
\begin{figure}[htbp]
\centering
  \includegraphics[height=3.5cm]{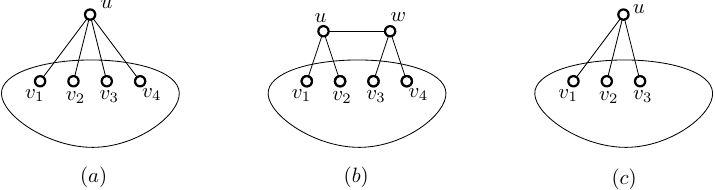}\\

  \caption{The last graphs in the sequence}
  \label{gi}

\end{figure}
\item[(b)] $G_1$ is as in \cref{gi}(b), with $v_1 \neq v_2$, $v_3 \neq v_4$, and $o_1$ is deleting $uw$.
\item[(c)] $G_1$ is as in \cref{gi}(c) and $o_1$ is to delete $uv_1$.
\item[(d)] $G_2$ is as in \cref{gi}(c) and $o_2$ is deletion of $uv_1$ (thus forming an edge $v_2v_3$), and $o_1$ is deletion of $v_2v_3$, so $G_2\setminus u \cong H$.
\end{itemize}

\end{lemma}
\begin{proof}
Let $G_k$ be the graph in the sequence which attains the $\min \{ k: |V(G_k)|>|V(H)| \}$, thus $V(G_{k-1})=V(H)=\{ v_1, v_2, \ldots v_{|H|} \}$. First, consider the case where $o_k$ is a split. Since this split reduces the number of vertices, it must be a split at a vertex $u$ of degree 4, see \cref{gi}$(a)$. Let $v_1v_2, v_3v_4$ be the edges resulting from splitting $u$. We claim that $G_{k-1}=H$, since if there was $k'<k$ so that $o_{k'}$ was
\begin{itemize}
\item splitting a 2-edge-path where both edges are present in $G_k$, or deleting an edge present in $G_k$, then it could have been done before $o_k$.
\item splitting a $v_1v_2v_i$ path, then we could have split $uv_2v_i$ instead.
\item splitting a 2-edge-path, with both edges $v_1v_2$ and $v_3v_4$, with, say, $v_2=v_3$, resulting from $o_k$, then we could have deleted one of $uv_2$ edges instead.
\item deleting one of the edges, say $v_1v_2$, created by $o_k$, then we could have deleted $uv_1$. (It also implies that $v_1\neq v_2$ and $v_3 \neq v_4$.)
\end{itemize}
Note that in all the cases above the alternative operation would result in another good sequence, with smaller $\min \{ k: |V(G_k)|>|V(H)| \}$, contradicting our choice of the good sequence. Therefore the claim is proved, thus $k=1$, and $(a)$ occurs.

Now, consider the case where $o_k$ is a deletion of an edge $uw$. Since this deletion reduces the number of vertices, at least one of its endpoints is of degree 3. If both $u$ and $w$ are of degree 3 (see \cref{gi}$(b)$), the same argument as above shows that $k=1$, and thus $(b)$ happens.

Otherwise, only $u$ is of degree 3, and $o_k$ is deleting $uv_1$, see \cref{gi}$(c)$. As before, it could be argued that there cannot be a $k'<k$ with $o_{k'}$ being splitting a 2-edge-path with both edges present in $G_k$, or deleting an edge present in $G_k$. Also, $o_{k'}$ cannot be splitting $v_2v_3v_i$, since we could have split $uv_3v_i$ before $o_k$, obtaining a good sequence with smaller $\min \{ k: |V(G_k)|>|V(H)| \}$. However, it could be that $o_{k'}$ is deleting $v_2v_3$. Thus, if $v_2v_3$ is not to be deleted, we have $k=1$ and $(c)$ happens; else, $k=2$ and $o_{k-1}$ would be deleting $v_2v_3$, i.e. $(d)$ occurs.
\end{proof}

Now, we use this lemma to establish a result on $K_5$-immersions discussed in \cref{intro} and {\color \col{repeated}} here.
\begin{cor}
\label{cor}
Suppose $G$ is $3$-edge-connected, and internally $4$-edge-connected, where $G\succ K_5$. Then
\begin{enumerate}
\item
if $|V(G)|\ge 6$ then $G\succeq K_{3,3}$, or $G\cong$ octahedron, where octahedron is the graph in \cref{Octahedron}.
\item
If $|V(G)|\ge 7$ then $G\succ K_{3,3}$.
\end{enumerate}
\end{cor}

\begin{proof}
Observe part $(2)$ is an immediate consequence of part $(1)$. We will then prove $(1)$.
Suppose $G\succ K_5$, and $|V(G)| >5$. By \cref{sp-thm2}, a good sequence from $G$ to $K_5$ exists. Thus, we can choose a good sequence
$$G=G_l, G_{l-1}, \ldots , G_2, G_1, G_0\cong K_5$$
such that $\min \{ k: |V(G_k)|>5 \}$ is as small as possible, and apply the previous lemma. It could be easily verified that if cases $(b)$ or $(c)$ of the previous lemma occur, then $G_1\succ K_{3,3}$, and if case $(d)$ happens, $G_2\succ K_{3,3}$, thus $G\succ K_{3,3}$.

So, suppose case $(a)$ of the previous lemma occurs. Again, it can easily be verified that if the two edges created by $o_1$ share an endpoint, then $G_1\succ K_{3,3}$, thus $G\succ K_{3,3}$. Otherwise, $K_{3, 3}$ is not immersed in $G_1$, as $G_1$ would be the octahedron, which, being planar, doesn't have $K_{3,3}$ as a subgraph. On the other hand, it has six vertices, all of degree 4, so an immersion of $K_{3,3}$ cannot be found doing splits either. 
\begin{figure}[htbp]
\centering
  \includegraphics[height=3cm]{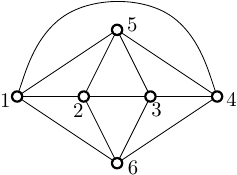}\\
  \caption{Octahedron}
  \label{Octahedron}
\end{figure}

Therefore, if $G\cong$ octahedron, $G\nsucceq K_{3,3}$. However, if $G$ properly immerses octahedron, then it immerses $K_{3,3}$ as well. To see that, note that the 6-vertex graphs from which octahedron is obtained after deletion of an edge or splitting a 2-edge path, all immerse $K_{3,3}$. On the other hand, if $|V(G)|>6$, we may again use \cref{lem11}, since octahedron itself is 4-edge-connected.

To reduce the number of graphs we examine, it now helps to notice that we only need to consider the case where a 4-vertex 7 gets split to create edges $\{23, 15\}$, or $\{ 23, 14\}$. It is because in all other cases, the graph obtained by splitting 2-paths 163, 264 would be one of the graphs we already looked at, all of which immerse $K_{3, 3}$.

If vertex 7 is split to create $\{23, 15\}$, an immersion of $K_{3,3}$ may be found after splitting 2-path 173. Also, if vertex 7 is split to create $\{23, 14\}$, then $K_{3, 3}$ lies as a subgraph in $G$.
\end{proof}
\section{Acknowledgment}
Both authors are supported in part by NSERC Discovery Grant RGPIN-06301.
\bibliographystyle{plain}
\bibliography{references}
\end{document}